\documentclass[reqno,a4paper,twoside,12pt]{amsart}
\usepackage{amsthm,amsmath,amsfonts}
\usepackage{graphicx}
\usepackage{hyperref}
\usepackage[latin1]{inputenc}
\theoremstyle{definition}
\newtheorem{definition}{Definition}

\theoremstyle{remark}

\theoremstyle{plain}
\newtheorem{lemma}[definition]{Lemma}

\newcommand{\set}[1]{\left\{{#1}\right\}}
\newcommand{\onemat}{\mathbb{I}}
\newcommand{\zeromat}{\mathbb{O}}

\newcommand{\dig}{\mathcal{DI}}
\newcommand{\pdig}{\mathcal{PDI}}

\newcommand\setsuchas[2]{\left\{\,{#1}\,\vrule\,{#2}\,\right\}}

\newcommand{\Q}{{\mathbb{Q}}}
\begin{document}
\title[Maximal spectral radius of a digraph]{The maximal spectral radius 
of a digraph with \((m+1)^2 - s\) edges}
\author{Jan Snellman}
\address{
Department of Mathematics, Stockholm University\\
SE-10691 Stockholm, Sweden}
\email{jans@math.su.se}

\date{March 10, 2003}
\begin{abstract}
  It is known that the  spectral radius of a digraph with \(k\)
  edges is \(\le \sqrt{k}\), and that this inequality is strict except
  when \(k\) is a perfect square. For \(k=m^2 + \ell\), 
  \(\ell\) fixed, \(m\)  large, Friedland
  showed that the optimal digraph is
  obtained from the complete digraph on \(m\) vertices by adding one
  extra vertex,  a corresponding loop,  and then 
  connecting it to the first \(\lfloor \ell/2\rfloor\)
  vertices by pairs of directed edges (for even \(\ell\) we add one
  extra edge to the new vertex). 

  Using a combinatorial reciprocity theorem by Gessel, and a
  classification by Backelin on the digraphs on \(s\) edges having
  a maximal number of walks of length two, we obtain the following
  result: for fixed \(0< s \neq 4\), \(k=(m+1)^2 - s\), \(m\) large,
  the maximal spectral  radius of a digraph with \(k\) edges 
  is obtained by the digraph which is constructed from the complete
  digraph on \(m+1\) vertices by removing the loop at the last vertex
  together with \(\lfloor s/2 \rfloor\)
  pairs of directed edges that connect to the last vertex (if \(s\) is
  even, remove an extra edge connecting to the last vertex).

  \textbf{Résumé}. \,\,
  On sait que le rayon spectral d'un graphe orienté avec \(k\) 
  arcs est \(\le \sqrt{k}\), et que cette inégalité est
  stricte sauf quand \(k\) est un carré parfait. Pour
  \(k=m^2 + \ell \), \(\ell\)  fixé, \(m\) grand, Friedland a
  prouvé que le graphe orienté optimal est obtenu à partir du graphe orienté
  complet à \(m\) sommets en ajoutant un sommet
  supplémentaire, une boucle correspondante, et en le reliant aux
  premiers \(\lfloor \ell/2\rfloor \)  sommets par des paires d'arcs.

  En utilisant un théorème
  combinatoire de réciprocité de Gessel, et une classification due à
  Backelin des graphes orientés à \(s\) arcs ayant un nombre
  maximal de chemins de  longueur deux, nous obtenons le
  résultat suivant: 
  pour \(s \neq 4\) fixé,   \(0< s\), 
  \(k=(m+1)^2 - s\), \(m\) grand, le rayon spectral maximal d'un
  graphe orienté à  
  \(k\)  arcs est obtenu pour le graphe orienté  construit à partir du
  graphe orienté complet à \(m+1\)  sommets en enlevant la boucle d'un
  sommet quelconque ainsi que les  \(\lfloor s/2 \rfloor\)  paires
  d'arcs reliée à ce  sommet.
  \end{abstract}

\thanks{Digraphs drawn by \textbf{dot} \cite{dot}.}  

\maketitle

\begin{section}{Introduction}
    By a \emph{digraph} we understand a finite directed graph with no
  multiple edges, but possibly loops.
  Let \(G(m,p,q)\) be the digraph on \(\set{1,\dots,m+1}\),
  where there is an edge from \(i\) to \(j\) if \(i,j \le m\) or if
  \(i \le p\) and \(j = m+1\) or if \(i=m+1\) and \(j \le q\). Let
  \(M(G(m,p,q)) = M(m,p,q)\) denote the adjacency matrix of
  \(G(m,p,q)\); it is a  0-1 matrix with \(m^2 + p + q\) ones.
  If \(\onemat_{a,b}\) denotes the \(a\times b\) matrix with
  all ones, and \(\zeromat_{a,b}\) the matrix with all zeroes,
  then
  \begin{equation}\label{eq:Mmpq}
    M(m,p,q)= \left[
    \begin{array}{c|c}
     \onemat_{m,m} & 
     \begin{smallmatrix}
       \onemat_{p,1}\\
       \zeromat_{m-p,1}
     \end{smallmatrix}
 \\ \hline
     \begin{smallmatrix}
       \onemat_{1,q} &  \zeromat_{1,m-q}
       \end{smallmatrix} & 0
    \end{array}\right].
  \end{equation}
  For \(0 < \ell < 2m+1\), we put \(M(m,\ell)=M(m,\lceil \ell/2
  \rceil, \ell - \lceil \ell/2 \rceil)\), i.e. the \(M(m,p,q)\)-matrix
  with \(p+q=\ell\),\(p\ge q\), \(p-q\) minimal. We denote the
  corresponding digraph with \(G(m,\ell)\).
  As an example, \(G(5,2,1)=G(5,3)\) is shown in
  figure~\ref{fig:G521}, and has adjacency matrix 
  \begin{equation}
    \label{eq:M521}
    M(5,2,1)=M(5,3) = 
    \begin{bmatrix}
      1 & 1 & 1 & 1 & 1 & 1 \\
      1 & 1 & 1 & 1 & 1 & 0 \\
      1 & 1 & 1 & 1 & 1 & 0 \\
      1 & 1 & 1 & 1 & 1 & 0 \\
      1 & 1 & 1 & 1 & 1 & 0 \\
      1 & 1 & 0 & 0 & 0 & 0 
    \end{bmatrix}
  \end{equation}
  \begin{figure}[htbp]
    \centering
  \includegraphics[scale=0.3, bb=18 18 594 516]{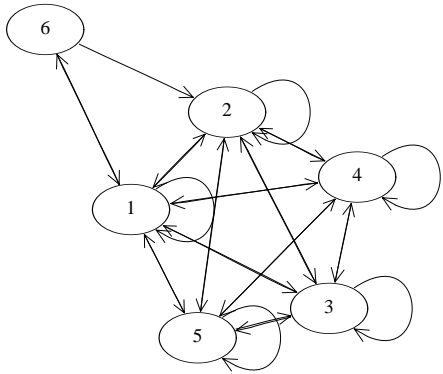}    
    \caption{ \(G(5,2,1)=G(5,3)\)}
    \label{fig:G521}
  \end{figure}

  Friedland \cite{Friedland:ME} showed that for a fixed \(\ell\) there
  is an \(L(\ell)\) 
  such that for \(m \ge L(\ell)\) the maximal spectral
  radius of a 0-1 matrix with \(m^2 + \ell\) ones is achieved by the
  matrix \(M(m,\ell)\).
  He conjectured that for any \(m,\ell\), \(0 < \ell < 2m +1\),
  the maximal spectral
  radius of a 0-1 matrix with \(m^2 + \ell\) ones is achieved by some
  \((m+1) \times (m+1)\) matrix.
  For \(m^2+\ell = (m+1)^2 - 4\),  he showed that the optimal
  matrix is \emph{not} \(M_{m,\ell}\) but rather 
  \begin{equation}\label{eq:M4}
    \left[
    \begin{array}{c|c}
      \onemat_{m-1,m-1} & \onemat_{m-1,2}\\ \hline
      \onemat_{2,m-1} & \zeromat_{2,2}
      \end{array}
    \right]
  \end{equation}
  It is reasonable to believe that this is the \emph{only exception} and
  that for other \(\ell\), \(M(m,\ell)\) is optimal.
  
  We show a weaker result: for a fixed \(s \neq 4\) there is an
  \(S(s)\) 
  such that for \(m \ge S(s)\) the maximal spectral
  radius of a digraph with \((m+1)^2 -s\)
  edges is achieved by the digraph \(G(m,2m + 1 - s)\).

  Our main tools are the following:
  \begin{enumerate}
  \item A combinatorial reciprocity theorem by Gessel \cite{Gessel}
    which asserts that for a digraph \(A\), the generating series
    \begin{equation}
      \label{eq:gHA}
      H_A(t) = \sum_{n=0}^\infty \chi_n(A)t^n,
    \end{equation}
    of \(\chi_n(A)\), the number of walks of length \(n-1\) in
    \(A\), is related to the series of  the complementary digraph by
    \(H_A(t) H_{\bar{A}}(-t)=1\),
  \item A classification by Backelin \cite{Backelin:NCM} of the
    digraphs of \(s\) edges with maximal number of walks of length 2.
  \end{enumerate}
  The proof runs as follows: Backelin's classification shows that for
  \(s >6\), \(m\) sufficiently large, the digraph with the following
  adjacency matrix has the 
  maximal number of walks of length 2 among digraphs with
  \(s\) edges and \(m+1\) vertices:
  \begin{equation}\label{eq:starmat}
 \left[
    \begin{array}{c|c}
     \zeromat_{m,m} & 
     \begin{smallmatrix}
       \zeromat_{m-s/2,1}\\
       \onemat_{s/2,1}
     \end{smallmatrix}
 \\ \hline
     \begin{smallmatrix}
       \zeromat_{1,m-s/2} &  \onemat_{1,m-s/2}
       \end{smallmatrix} & 1
    \end{array}\right].
  \end{equation}
  The  generating series for walks in that graph is \[1+(m+1)t + st^2
  + ct^3 + O(t^4),\] 
  so the generating series for the complementary graph, which has
  adjacency matrix \(M(m,2m+1-s)\), is 
  \begin{displaymath}
    \frac{1}{1-(m+1)t + st^2 - ct^3 + O(t^4)}.
  \end{displaymath}
  For any other digraph with \(s\) edges and \(m+1\) vertices, we have
  that the generating series is 
  \begin{displaymath}
    \frac{1}{1-(m+1)t + st^2 - dt^3 + O(t^4)},
  \end{displaymath}
  with \(d<c\). A perturbation analysis yields that  pole of smallest
  modulus is located at \(m^{-1} + sm^{-3} -cm^{-4}\) in the first
  case and at 
  \(m^{-1} + sm^{-3} -dm^{-4}\) in the second, so the first series has
  smaller radius of convergence, hence faster growth of the
  coefficients. Consequently, the first graph has the larger spectral
  radius. 
\end{section}

\begin{section}{The proof}
  For any digraph \(A\), let 
  \begin{equation}
    \label{eq:HA}
    H_A(t) = \sum_{n=0}^\infty \chi_n(A)t^n,
  \end{equation}
  where \(\chi_n(A)\) denotes  the number of walks of length \(n\) in
  \(A\), so that \(\chi_0(A)=1\), \(\chi_1(A)=\) the number of
  vertices in \(A\), and \(\chi_1(A)=\) the number of edges in \(A\).
  Let 
  \begin{equation}
    \label{eq:RA}
R(A) = \frac{1}{\lim_{n \to \infty} \sqrt[n]{\chi_n(A)}}
  \end{equation}
 be the  radius of convergence of \(A\), and let \(\rho(A)=1/R(A)\).
 If the adjacency matrix \(M(A)\) of \(A\) is irreducible then
 \(\rho(A)>0\) is 
 the largest eigenvalue of \(M(A)\).

 For \((m+1)^2-s > m^2+1\), let \(\dig(m,s)\) denote the finite set of
 digraphs on 
 \(\set{1,\dots,m+1}\) 
 having precisely \((m+1)^2 - s\) edges. 
 Let \(\pdig(m,s) \subset \dig(m,s)\) denote the subset consisting of
 those digraphs whose 0-1 adjacency can be regarded as the Young
 diagram of a numerical partition of \((m+1)^2 - s\); in other words, 
 the rows and columns of the adjacency matrix should be weakly
 decreasing. 
 Then \(\pdig(m,s)\) is finite, and the cardinality does not depend on
 \(m\) as long as \(m\) is sufficiently large.
 Furthermore every digraph in \(\pdig(m,s)\) is connected in the
 directed sense, i.e. there is a directed walk between any two vertices.
 digraphs. Hence the
 adjacency matrix of an element in \(\pdig(m,s)\) is by definition
 \emph{irreducible}.
 Furthermore, by a result of Schwarz \cite{Schwarz:Re}, 
 \begin{equation}
   \label{eq:sc}
   \max \setsuchas{\rho(A)}{A \in \dig} = \max \setsuchas{\rho(A)}{A
     \in \pdig(m,s)}. 
 \end{equation}

 Let \(\bar{A}\) denote the complementary graph of \(A\),
 i.e. the digraph on \(\set{1,\dots,m+1}\) which has an edge  \(i \to
 j\) iff there isn't an edge \(i \to j\) in \(A\). 
 Then the following relation hold (see \cite{SyzygiesAndWalks} and
 also \cite{Gessel} and \cite{EnPair}): 
 \begin{equation}
   \label{eq:duality}
   H_A(t)H_{\bar{A}}(-t) = 1
 \end{equation}
 If \(A \in \pdig(m,s)\), then \(\bar{A}\) is a digraph on \(m+1\)
 vertices 
 with \(s\) edges. 
 We have that 
 \begin{equation}
   \label{eq:bara}
   H_{\bar{A}}(t) = 1 + (m+1)t + st^2 + ct^3 + O(t^4)
 \end{equation}
 where \(c\) is the number of walks in \(\bar{A}\) of length \(2\).

 \begin{subsection}{The case \(s>6\)}
 Suppose first that \(s>6\). Backelin \cite{Backelin:NCM} showed that
 among all digraphs with \(s\) edges, the so-called \emph{saturated
   stars} have the maximal number of walks of length 2. By a
 saturated star with \(s=2k-1\) edges is meant the digraph with a
 edges 
 \((1,i)\) and \((i,1)\) for \(1 \le i \le k\); for \(s=2k\) we add
 the 
 edge \((1,k+1)\). 
 So the saturated stars with 9 and 10 edges looks as in
 figure~\ref{fig:satstar}
 \begin{figure}[htbp]
   \centering
 \begin{tabular}{cc}
\includegraphics[scale=0.5]{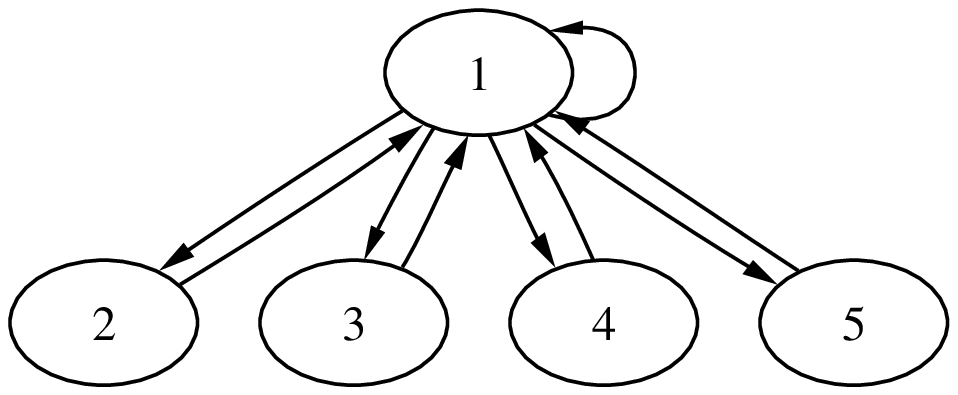} &
\includegraphics[scale=0.5]{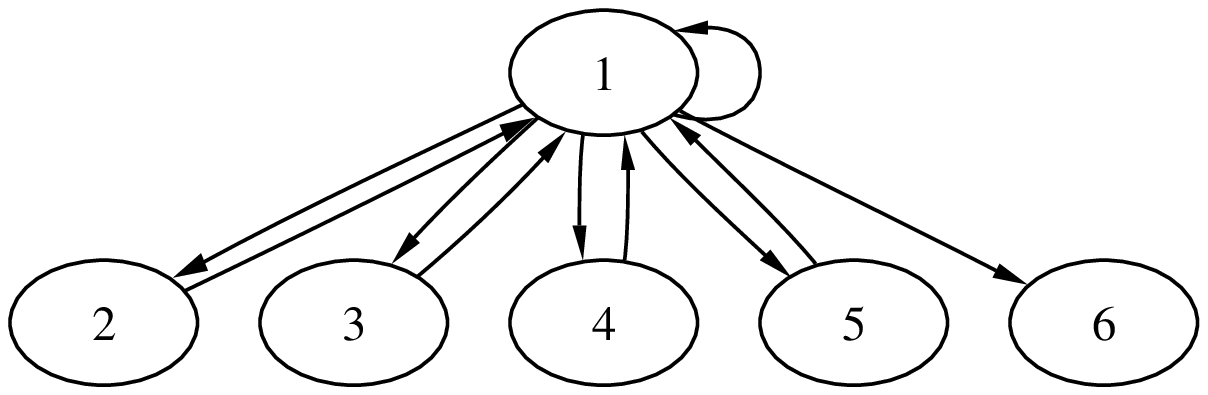}
 \end{tabular}
   \caption{Saturated star digraphs with 9 and 10 edges}
   \label{fig:satstar}
 \end{figure}

Note that if \(R\) is a saturated star, then the graph on \(m+1\)
vertices, which has an edge \(i\to j\) iff \((m+2-i) \to (m+2-j)\) is not an
edge in \(R\), is of the form \(G(m,s)\). For instance, the digraphs
above have  adjacency matrices
 \begin{equation}
   \label{eq:adjstar}
   \begin{bmatrix}
     1 & 1 & 1 & 1 & 1\\
     1 & 0 & 0 & 0 & 0\\
     1 & 0 & 0 & 0 & 0\\
     1 & 0 & 0 & 0 & 0\\
     1 & 0 & 0 & 0 & 0
   \end{bmatrix}
   \quad \text{ and }
      \begin{bmatrix}
     1 & 1 & 1 & 1 & 1 & 1\\
     1 & 0 & 0 & 0 & 0 & 0\\
     1 & 0 & 0 & 0 & 0 & 0\\
     1 & 0 & 0 & 0 & 0 & 0\\
     1 & 0 & 0 & 0 & 0 & 0\\
     0 & 0 & 0 & 0 & 0 & 0
   \end{bmatrix}
 \end{equation}
and if we take \(m=7\) we get the following adjacency matrices for the
relabeled complementary graphs:
\begin{equation}
  \label{eq:adjcomp}
   \begin{bmatrix}
     1 & 1 & 1 & 1 & 1 &1 & 1 & 1 \\
     1 & 1 & 1 & 1 & 1 &1 & 1 & 1 \\
     1 & 1 & 1 & 1 & 1 &1 & 1 & 1 \\
     1 & 1 & 1 & 1 & 1 &1 & 1 & 0 \\
     1 & 1 & 1 & 1 & 1 &1 & 1 & 0 \\
     1 & 1 & 1 & 1 & 1 &1 & 1 & 0 \\
     1 & 1 & 1 & 1 & 1 &1 & 1 & 0 \\
     1 & 1 & 1 & 0 & 0 &0 & 0 & 0 
   \end{bmatrix}
   \quad \text{ and }
      \begin{bmatrix}
     1 & 1 & 1 & 1 & 1 &1 & 1 & 1 \\
     1 & 1 & 1 & 1 & 1 &1 & 1 & 1 \\
     1 & 1 & 1 & 1 & 1 &1 & 1 & 1 \\
     1 & 1 & 1 & 1 & 1 &1 & 1 & 0 \\
     1 & 1 & 1 & 1 & 1 &1 & 1 & 0 \\
     1 & 1 & 1 & 1 & 1 &1 & 1 & 0 \\
     1 & 1 & 1 & 1 & 1 &1 & 1 & 0 \\
     1 & 1 & 0 & 0 & 0 &0 & 0 & 0 
   \end{bmatrix}
\end{equation}
 
Now suppose that \(R\) is a saturated star, and that \(A\) is the
digraph on \(m+1\) vertices obtained from \(R\) as above. So 
\(\bar{A}\) and \(R\) differ only in that \(\bar{A}\) has some
isolated vertices.
Let \(R'\) be a
a different digraph with \(s\) edges. Let \(\bar{B}\) be the digraph
obtained by adjoining isolated vertices, so that the total number of
vertices becomes \(m+1\). Let \(B = \bar{\bar{B}}\). Then 
 \begin{equation}
   \label{eq:barb}
   H_{\bar{B}}(t) = 1 + (m+1)t + st^2 + dt^3 + O(t^4)
 \end{equation}
 where \(d\) is the number of walks in \(R'\) of length \(2\).
 By Backelin's result,
 \(c>d\). Hence we have that 
 \begin{equation}
   \label{eq:diff}
   \begin{split}
     H_{A}(t) - H_{B}(t) &=
     \frac{1}{1 - (m+1)t + st^2 - ct^3 + O(t^4)}
      \\
     & \qquad \qquad - \frac{1}{1 - (m+1)t + st^2 - dt^3 + O(t^4)} \\
     &= (c-d)t^3 + O(t^4)
   \end{split}
 \end{equation}
so we have at once that \(A\) has strictly more walks
of length 2 than \(B\) has. 
By induction, we can show that 
\begin{lemma}\label{lemma:expo}
The exponent of
\(t^i\) in \eqref{eq:diff} is a polynomial in \(m\) of degree \(i-3\),
with leading coefficient \((i-2)(c-d)\).  
\end{lemma}

 Thus, for any \(j\), by
taking \(m\) sufficiently large, we can achieve that the coefficients 
of \(t^i\), \(i \le j\), in \eqref{eq:diff} are all positive.
Recall that \(\pdig(m,s)\) is finite. Hence, for any \(j\), if we take
\(m\) sufficiently large, then \(A\) has the maximal number of walks
of length \(i\) among the \(B \in \pdig(m,s)\).

In \cite{Stanley:En1} it is shown that if \(\mathfrak{G}\) is a
digraph with adjacency matrix \(M\), 
then \(H_{\mathfrak{G}}(t) = P(t)/Q(t)\), where \(Q(t)=\det(I - tM)\),
and \(P(t)\) is a polynomial of smaller degree then \(Q\).
Hence \(H_{\bar{A}}(t), H_{\bar{B}}(t), H_{A}(t), H_{B}(t)\) are all
rational functions.
Let \(r=r(m,c)\) be the pole of \(H_A(t)\) that is closest to
origin. Then \(1/r\) is the eigenvalue of  \(M\) of largest modulus, so from the
Perron-Frobenius theorem  it follows that if \(\mathfrak{G}\) is
connected in the directed sense then \(r\) is a positive real number. Writing 
\begin{equation}
  H_A(t) = \frac{1}{1 -(m+1)t + st^2- ct^3 + t^4\frac{b_0 + b_1t +
      \cdots + b_{N_1}t^{N_1}}{1+a_1t + 
    a_2t^2 + \cdots + a_{N_2}t^{N_2}}} 
\end{equation}
we have that  \(r(m,c)\) is the smallest real root of 
\begin{multline}\label{eq:sr}
(1+a_1t +  a_2t^2 + \cdots + a_{N_2}t^{N_2})(1-(m+1)t+st^2-ct^3) 
\\ +
t^4(b_0 
+ b_1t + \cdots + b_{N_1}t^{N_1})=0. 
\end{multline}
We are interested in the asymptotic behavior of \(r\) as \(m\to
\infty\).
As we will show below, we can
expand \(r(m,c)\) as a Laurent  series in \(m\) as
\begin{equation}
  \label{eq:expansion}
  r(m,c)= m^{-1} + sm^{-3} -cm^{-4} +
  (2s^2+b_0)m^{-5} + 
(b_1-a_1b_0 -  5cs)m^{-6} + O(m^{-7})
\end{equation}
Making a similar analysis for the pole \(r^*(m,d)\) of
\(H_B(t)\) that has the smallest modulus, we 
get \[r(m,c) - r^*(m,d) = -(c-d)m^{-4} +
O(m^{-5}),\] so 
\[r(m,c) < r^*(m,d).\] Thus, for large \(m\), \(H_A(t)\) has
strictly smaller radius of convergence than
\(H_B(t)\). When we combine this with Lemma~\ref{lemma:expo}
we see that 
by taking \(m\) sufficiently large, we can achieve that 
  \(H_A(t) \gg H_B(t)\),
i.e. that \emph{all} coefficients of \(H_A(t)\) are \(\ge\)
than the corresponding coefficients of \(H_B(t)\). In fact,
the inequality is strict for exponents \(>2\).

\begin{subsubsection}{Perturbation analysis of the positive root}
We now show how to derive the expansion   \eqref{eq:expansion}.
Replacing \(m+1\) by \(1/\varepsilon\) in \eqref{eq:sr}, and clearing
denominators, we get 
\begin{multline}
  \label{eq:cd}
  (1+a_1t +  a_2t^2 + \cdots + a_{N_2}t^{N_2})(\varepsilon-t+\varepsilon
  st^2-\varepsilon ct^3) \\ + \varepsilon t^4(b_0
+ b_1t + \cdots + b_{N_1}t^{N_1}) =0
\end{multline}
The unperturbed equation is 
\begin{equation}
  \label{eq:cd2}
   (1+a_1t +  a_2t^2 + \cdots + a_{N_2}t^{N_2})(-t) = 0,
\end{equation}
which has a root at \(t=0\). We 
introduce the scaling
\(t=\varepsilon T\) and get 
\begin{multline}
  \label{eq:cd3}
  (1+a_1\varepsilon T +  a_2\varepsilon^2T^2 + \cdots + a_{N_2}
  \varepsilon^{N_2}T^{N_2})(\varepsilon-\varepsilon T+\varepsilon^3 
  sT^2-\varepsilon^4 cT^3) \\+ \varepsilon^5 T^4(b_0
+ b_1\varepsilon T + \cdots + b_{N_1}\varepsilon^{N_1}T^{N_1}) =0
\end{multline}
hence
\begin{multline}
  \label{eq:cd4}
  (1+a_1\varepsilon T +  a_2\varepsilon^2T^2 + \cdots + a_{N_2}
  \varepsilon^{N_2}T^{N_2})(1- T+\varepsilon^2 
  sT^2-\varepsilon^3 cT^3) \\+ \varepsilon^4 T^4(b_0
+ b_1\varepsilon T + \cdots + b_{N_1}\varepsilon^{N_1}T^{N_1}) =0
\end{multline}
The unperturbed equation is now \(1-T=0\). Hence, \(T=O(1)\) and
\(T^{-1} =O(1)\), so this is the correct scaling.
We make the substitution \(Y=T-1\) and get
\begin{multline}
  \label{eq:cd5}
  \bigl[1+a_1\varepsilon (Y+1) +  a_2\varepsilon^2(Y+1)^2 + \cdots + a_{N_2}
  \varepsilon^{N_2}(Y+1)^{N_2}\bigr]\bigl[-Y+\varepsilon^2 
  s(Y+1)^2- \\
\varepsilon^3 c(Y+1)^3\bigr] + \varepsilon^4 (Y+1)^4(b_0
+ b_1\varepsilon (Y+1) + \cdots + b_{N_1}\varepsilon^{N_1}(Y+1)^{N_1}) =0
\end{multline}
It is now clear that \(Y\) can be expanded in powers of \(\varepsilon\),
so we make the \emph{Ansatz} 
\begin{equation}
  \label{eq:Ansatz}
  Y=\sum_{i=1}^\infty w_i\varepsilon^i.
\end{equation}
Collecting the coefficients of the powers of \(\varepsilon\) in 
\eqref{eq:cd5} we get 
\begin{eqnarray}
  \label{eq:cd6}
  1: & 0 \\
  \varepsilon: & -w_1 \\
  \varepsilon^2: &  s-a_1w_1-w_2  \\
  \varepsilon^3: &  2sw_1+ a_1s-c-a_1w_1^2- a_1w_2-a_2w_1-w-3  \\
  \varepsilon^4: &   2sw_2 + sw_1^2 -3cw_1- a_1c +a_2s \\
&  +3a_1sw_1 +b_0
-2a_1w_1w_2
-a_1w_3 -2a_2w_1^2 -a_2w_2 -a_3w_1 -w_4 
\end{eqnarray}
These should be zero, which allows us to solve for the \(w_i\)'s,
obtaining
\begin{equation}
  \label{eq:b1}
  w_1 = 0, \quad w_2 = s, \quad w_3 = -c, \quad w_4 = 2s^2 +b_0
\end{equation}
So 
\begin{equation}
  \label{eq:b2}
  Y= s\varepsilon^2 -c\varepsilon^3 + (2s^2 +b_0)\varepsilon^4 +
  O(\varepsilon^5), 
\end{equation}
hence
\begin{equation}
  \label{eq:b3}
  T= 1 +s\varepsilon^2 -c\varepsilon^3 + (2s^2 +b_0)\varepsilon^4
  +O(\varepsilon^5), 
\end{equation}
hence
\begin{equation}
  \label{eq:b4}
 t= \varepsilon +s\varepsilon^3 -c\varepsilon^4 + (2s^2
 +b_0)\varepsilon^5 
 +O(\varepsilon^6). 
\end{equation}
\end{subsubsection}

This concludes the proof for the case  \(s >6\).
   
 \end{subsection}

\begin{subsection}{The exceptional cases}
It remains to take care of the case \(s \le 6\). 
 Backelin's classification says that if  \(s\in \set{1,3,5}\) then the
 saturated stars are optimal. Hence, it remains to check \(s=2\),
 \(s=4\), \(s=6\).

\begin{subsubsection}{\(s=2\)}
For \(s=2\) there are two non-isomorphic graphs \(R\),
namely 
\begin{center}
\includegraphics[scale=0.5]{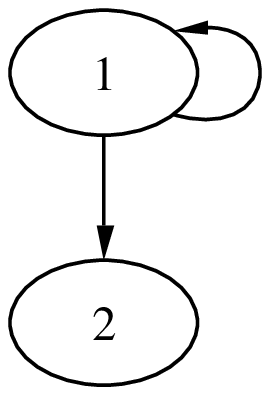} 
\includegraphics[scale=0.5]{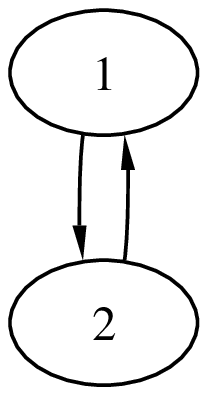}   
\end{center}
corresponding
to matrices 
\(M_1 = \bigl[\begin{smallmatrix}1 &1\\0 &0\end{smallmatrix}\bigr]\)
and  
\(M_2 = \bigl[\begin{smallmatrix}0 &1\\1
  &0\end{smallmatrix}\bigr]\). Recall that when \(M_1\) or \(M_2\) is
subtracted from the bottom right corner of the \((m+1)\times(m+1)\)
matrix of all ones, the result should be weakly decreasing in rows and
columns. We call the resulting matrices 
\begin{equation}
  \label{eq:A1}
  A_1 = 
  \begin{bmatrix}
    1 & \cdots & 1 &1 & 1\\
    \vdots & \vdots & \vdots & \vdots & \vdots\\
    1 & \cdots & 1 &1 & 1\\
    1 & \cdots & 1& 0 & 0
  \end{bmatrix}
  \qquad \qquad
  A_2 = 
  \begin{bmatrix}
    1 & \cdots & 1 &1 & 1\\
    \vdots & \vdots & \vdots & \vdots & \vdots\\
    1 & \cdots & 1 &1 & 0\\
    1 & \cdots & 1& 0 & 1
  \end{bmatrix}
\end{equation}
The matrix \(A_2\) is not weakly decreasing in the last row, so it is
not really necessary to continue with the calculations, but we proceed
anyway in order to demonstrate how this is done.
We have that
\begin{equation}
  \label{eq:M12}
  H_{M_1}(t) = H_{M_2}(t)=\frac{1+t}{1-t},
\end{equation}
so that
\begin{equation}
  \label{eq:A12}
  H_{A_1}(t) = H_{A_2}(t)=\frac{1}{\frac{1-t}{1+t} -m +1}
  = \frac{1+t}{1-tm -(m-1)t^2}
\end{equation}
The smallest positive root of the denominator is 
\begin{equation}
  {\frac {-m+\sqrt {{m}^{2}+4\,m-4}}{2(m-1)}},
\end{equation}
so the spectral radius is 
\begin{equation}
  \frac {2(m-1)}{-m+\sqrt {{m}^{2}+4\,m-4}}.
\end{equation}
  
\end{subsubsection}

\begin{subsubsection}{\(s=4\)}
From Backelin's classification we have that for \(s=4\),
the digraph
\begin{center}
  \includegraphics[scale=0.5]{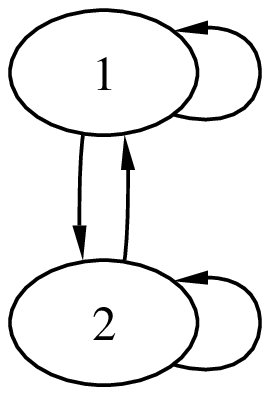}
\end{center}
with adjacency matrix
\(\bigl[\begin{smallmatrix}1 &1\\1 &1\end{smallmatrix}\bigr]\) has the
most walks of length 2. This is in accordance with \cite{Friedland:ME}
where it is shown that the matrix \eqref{eq:M4}
has maximal spectral radius among 0-1 matrices with \((m+1)^2 - 4\)
ones.
\end{subsubsection}

\begin{subsubsection}{\(s=6\)}

The remaining exceptional case in Backelin's classification is for
\(s=6\). 
\begin{figure}[htbp]
  \centering
  \begin{tabular}{cc}
    \includegraphics[scale=0.5]{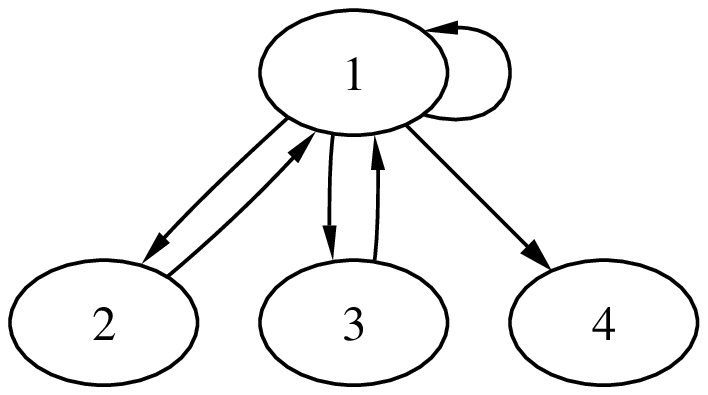} &
    \includegraphics[scale=0.5]{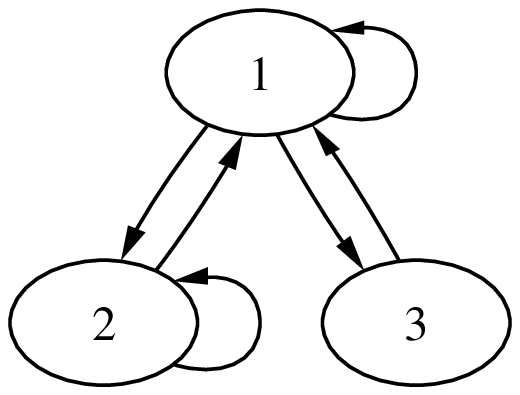} 
  \end{tabular}
  \caption{Digraphs with 6 edges having maximal number of walks of
    length 2} 
  \label{fig:s6}
\end{figure}
Then 
the  digraphs in figure~\ref{fig:s6}, with adjacency matrices
\begin{equation}
  \label{eq:M6}
    M_1 = 
    \begin{bmatrix}
      1 & 1 & 1 & 1 \\
      1 & 0 & 0 & 0 \\
      1 & 0 & 0 & 0 \\
      0 & 0 & 0 & 0
    \end{bmatrix}
    , \qquad
    M_2 = 
    \begin{bmatrix}
      1 & 1 & 1  \\
      1 & 1 & 0  \\
      1 & 0 & 0 
    \end{bmatrix}
\end{equation}
both have 14 of walks of length 2. Note that the first digraph is a
saturated star.
The relabeled complemented matrices are 
\begin{equation}
  \label{eq:A126}
  A_1 = 
  \begin{bmatrix}
    1 & \cdots & 1 &1 & 1 & 1 & 1\\
    \vdots & \vdots & \vdots & \vdots & \vdots & \vdots & \vdots\\
    1 & \cdots & 1 &1 & 1 & 1 & 1\\
    1 & \cdots & 1 &1 & 1 & 1 & 0\\
    1 & \cdots & 1 &1 & 1 & 1 & 0\\
    1 & \cdots & 1& 0 & 0 &0 & 0   
  \end{bmatrix}
  \qquad
  A_2 =
  \begin{bmatrix}
    1 & \cdots & 1 &1 & 1 & 1 \\
    \vdots & \vdots & \vdots  & \vdots & \vdots & \vdots\\
    1 & \cdots & 1 &1 & 1 & 1 \\
    1 & \cdots & 1 &1 & 1 & 0 \\
    1 & \cdots & 1 &1 & 0 & 0 \\
    1 & \cdots & 1& 0 & 0 &0    
  \end{bmatrix}
\end{equation}

The generating series for the complemented relabeled digraphs are 
\begin{eqnarray}
\label{eq:hser6}
    H_{A_1}(t) &= \frac{1 + t -2t^2}{ 1- mt -mt^2 +3t^2  + 2mt^3-6t^3}
    &=: \frac{p_1(t)}{p_2(t)}
    \\ 
    H_{A_2}(t) &= \frac{1+2t-t^2-t^3}{1 + t -mt +3t^2 -2mt^2  -2t^3
      + 
    mt^3 -2t^4 + mt^4}
  &=: \frac{q_1(t)}{q_2(t)}
\end{eqnarray}
Regarding the positive root of \(p_2(t)\) as a function of \(m\), and
expanding that function as a power series round infinity, we get 
\begin{displaymath}
  r_1(m) := m^{-1} - m^{-2} + 7m^{-3} -33m^{-4} + 191m^{-5} +
  O(m^{-6}), 
\end{displaymath}
whereas the expansion of the positive root of \(q_2(t)\) is
\begin{displaymath}
  r_2(m):=m^{-1} - m^{-2} + 7m^{-3} -33m^{-4} + 196m^{-5} + O(m^{-6}).
\end{displaymath}
The first root is therefore slightly smaller for large \(m\); the
difference is \emph{miniscule}, but \emph{vive la diff\'erence}! 
In fact, since \(r_1(m)-r_2(m)<0\) for large \(m\), and
\(r_1(4)-r_2(4) \approx -0.003\), it will suffice to show that
\(p_2(t)=q_2(t)=0\) has no solution to demonstrate that
\(r_1(m)-r_2(m)<0\) for \(m \ge 4\), as shown in
Figure~\ref{fig:diffe}.
 Using Macaulay 2 \cite{MACAULAY2}
we can verify that \(\set{1}\) is a Gr\"obner bases for the ideal
generated by  \(\set{p_2(t), q_2(t)}\) in \(\Q(m)[t]\).
Hence  \(p_2(t), q_2(t)\) are
co-prime in \(\Q(m)[t]\),  so they can not have a common zero.
\begin{figure}[htbp]
  \centering
  \caption{Difference \(r_1(m)-r_2(m)\) for  \(4 \le m \le 10\).}
  \includegraphics[scale=0.4, angle=270, bb= 70 70 530 670]{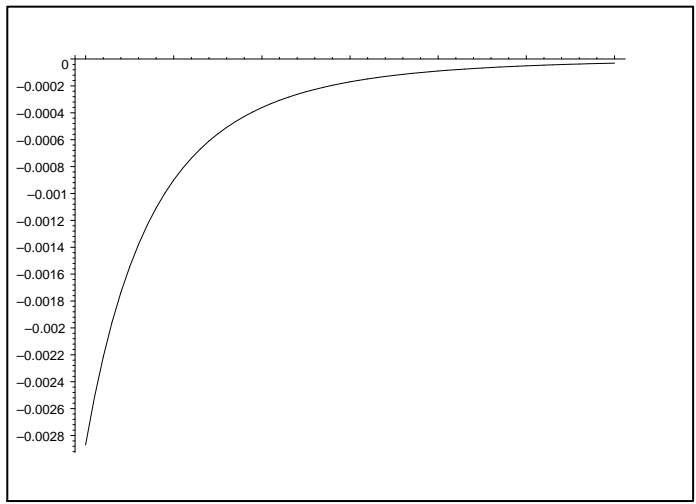}     
  \label{fig:diffe}
\end{figure}
\end{subsubsection}
\end{subsection}

\end{section}
\sloppy
\bibliographystyle{amsalpha}
\bibliography{journals,entropy,articles}
\end{document}